\documentclass[dvips]{siamltex}

\usepackage[utf8]{inputenc} % Use smart codification
\usepackage{psfrag}
\usepackage{epsfig}
\usepackage{bf1}
\usepackage{amsmath}
\usepackage{amssymb}
\usepackage{amsfonts}
\usepackage[section]{algorithm}
\usepackage{algpseudocode} % https://tex.stackexchange.com/questions/229355/algorithm-algorithmic-algorithmicx-algorithm2e-algpseudocode-confused
\usepackage[table]{xcolor}
\usepackage{tabularx} % Better control of columns size
\usepackage{booktabs} % Fancy tables from the guide https://www.inf.ethz.ch/personal/markusp/teaching/guides/guide-tables.pdf
\usepackage{multirow}
\usepackage{color}
\usepackage[english]{babel}
\renewcommand\arraystretch{1}

\usepackage{moreverb,url}
\usepackage{verbatimbox}
\usepackage{listings}
\usepackage{comment}
\title{A GPU-accelerated Adaptive FSAI preconditioner For Massively Parallel Simulations}
\author{Giovanni Isotton \footnotemark[1] \and Carlo Janna\footnotemark[1] \and
        Massimo Bernaschi\footnotemark[2]}

\begin{document}
\maketitle

%\runninghead{A GPU-accelerated Adaptive FSAI preconditioner for Massively Parallel Simulations}

\renewcommand{\thefootnote}{\fnsymbol{footnote}}
\footnotetext[1] {M3E S.r.l., via Giambellino 7, 35129 Paova, Italy,
{\tt e-mail} g.isotton@m3eweb.it,\,c.janna@m3eweb.it}
\footnotetext[2] {Institute for Applied Computing, CNR, 00185 Rome, Italy,
{\tt e-mail} massimo.bernaschi@cnr.it}
\renewcommand{\thefootnote}{\arabic{footnote}}

\begin{abstract}
The solution of linear systems of equations is a central task in a number of scientific
and engineering applications. In many cases the solution of linear systems may take most
of the simulation time thus representing a major bottleneck in the further development of
scientific and technical software.
For large scale simulations, nowadays accounting for several millions
or even billions of unknowns, it is quite common to resort to
preconditioned iterative solvers for exploiting their low memory
requirements and, at least potential, parallelism. Approximate
inverses have been shown to be robust and effective preconditioners in
various contexts.  In this work, we show how adaptive FSAI, an
approximate inverse characterized by a very high degree of
parallelism, can be successfully implemented on a distributed memory
computer equipped with GPU accelerators. Taking advantage of GPUs in
adaptive FSAI set-up is not a trivial task, nevertheless we show
through an extensive numerical experimentation how the proposed
approach outperforms more traditional preconditioners and results in a
close-to-ideal behaviour in challenging linear algebra problems.
\end{abstract}


\newcommand{\bb}{\mbox{\boldmath $ b$}}
\newcommand{\be}{\mbox{\boldmath $ e$}}
\newcommand{\bg}{\mbox{\boldmath $ g$}}
\newcommand{\bh}{\mbox{\boldmath $ h$}}
\newcommand{\br}{\mathbf{r}}
\newcommand{\bv}{\mathbf{v}}
\newcommand{\bw}{\mbox{\boldmath $ w$}}
\newcommand{\bx}{\mbox{\boldmath $ x$}}
\newcommand{\by}{\mbox{\boldmath $ y$}}
\newcommand{\bz}{\mbox{\boldmath $ z$}}
\newcommand{\bA}{\mbox{\boldmath $ A$}}
\newcommand{\bvartheta}{\mbox{\boldmath $ \vartheta$}}
\newcommand{\bzero}{\mbox{\boldmath $0$}}
\newcommand{\bnabla}{\mbox{\boldmath $\nabla$}}
\newcommand{\dt}{\mathrm{det}}
\newcommand{\tr}{\mathrm{tr}}
\newcommand{\diagb}{\mathrm{diag}_{\mathrm{B}}}
\newcommand{\ddiag}{\mathrm{diag}}
\newcommand{\daij}{ \frac{\displaystyle \partial}{\displaystyle \partial [A]_{ij}} }
\newcommand{\dfij}{ \frac{\displaystyle \partial}{\displaystyle \partial [F]_{ij}} }
\newcommand{\dFAFTa}{ \frac{\displaystyle \partial [FAF^T]_{ii}}{\displaystyle \partial
[F]_{ij}} }
\newcommand{\dFAFTb}{ \frac{\displaystyle \partial [FAF^T]_{ii}}{\displaystyle \partial
\mathbf{f}_{\PP_i} } }
\newcommand{\bndA}{\displaystyle \Psi_F(A)}
\newcommand{\blk}{\displaystyle \mathrm{Low}}
\newcommand{\I}{\mathcal{I}}
\newcommand{\J}{\mathcal{J}}
\newcommand{\JJ}{\hat{\mathcal{J}}}
\newcommand{\PP}{\mathcal{P}}
\newcommand{\TPP}{\widetilde{\mathcal{P}}}
\newcommand{\RR}{\mathcal{R}}
\newcommand{\CC}{\mathcal{C}}
\newcommand{\tbg}{\widetilde{\bg}}
\newcommand{\hbg}{\widehat{\bg}}
\newcommand{\Grc}{\widehat{G}[\RR,\CC]}
\newcommand{\Grct}{\widehat{G}[\RR,\CC]^T}
\newcommand{\Acc}{A[\CC,\CC]}
\newcommand{\Acr}{A[\CC,\RR]}
\newcommand{\Acrt}{A[\CC,\RR]^T}
\newcommand{\Acj}{A[\CC,\overline{j}]}
\newcommand{\Acjt}{A[\CC,\overline{j}]^T}
\newcommand{\Ajr}{A[\overline{j},\RR]}
\newcommand{\Ajrt}{A[\overline{j},\RR]^T}
\newcommand{\Arr}{A[\RR,\RR]}
\newcommand{\Ajj}{A[\overline{j},\overline{j}]}
\newcommand{\grc}{\widehat{G}_{rc}}
\newcommand{\acc}{A_{cc}}
\newcommand{\acr}{A_{cr}}
\newcommand{\arr}{A_{rr}}
\newcommand{\tg}{\widetilde{G}}
\newcommand{\hg}{\widehat{G}}
\newcommand{\oj}{\overline{j}}
\newcommand{\oi}{\overline{i}}
\newcommand{\aj}{\ba_{\oj}}
\newcommand{\ai}{\ba_{\oi}}
\newcommand{\gi}{\hbg_{\oi}}
\newcommand{\ajj}{a_{\oj\oj}}
\newcommand{\aii}{a_{\oi\oi}}
\newcommand{\bii}{b_{\oi\oi}}
\newcommand{\aji}{a_{\oj\oi}}

\newcommand{\rt}[1]{\textcolor{red}{#1}}
\newcommand{\bt}[1]{\textcolor{blue}{#1}}
\newcommand{\gt}[1]{\textcolor{green}{#1}}

% Carlo
\newcommand{\bbv}{\overline{\bv}}
\newcommand{\Wid}{W_{\mbox{\scriptsize ideal}}}
\newcommand{\bdelta}{\boldsymbol \delta}
\newcommand{\amgname}{aSP-AMG}

% Victor
\newcommand{\ra}[1]{\renewcommand{\arraystretch}{#1}}
\newcommand{\bigo}[1]{\mathcal{O}\left(#1\right)}
\renewcommand{\vec}[1]{\mathbf{#1}}
\providecommand{\vecT}[1]{\mathbf{#1}^{T}}
\providecommand{\mat}[1]{#1}
\providecommand{\matTilde}[1]{\tilde{#1}}
\providecommand{\matT}[1]{{#1}^{T}}
\providecommand{\matI}[1]{{#1}^{-1}}
\providecommand{\matIT}[1]{{#1}^{-T}}
\providecommand{\rap}[3]{\mat{#1} \mat{#2} \mat{#3}}
\providecommand{\rapI}[3]{\left(\mat{#1} \mat{#2} \mat{#3} \right)^{-1}}
\providecommand{\ptap}[2]{\mat{#1}^{T} \mat{#2} \mat{#1}}
\providecommand{\ptapI}[2]{\left(\mat{#1}^{T} \mat{#2} \mat{#1} \right)^{-1}}
\providecommand{\fronorm}[1]{\left|| #1 |\right|_{F}}
\providecommand{\twonorm}[1]{\left|| #1 |\right|_{2}}
\providecommand{\gradsym}[1]{\left( \nabla \vec{#1} + \nabla \vecT{#1} \right)}
\providecommand{\divB}[1]{\nabla \cdot \left[ #1 \right]}
\newcommand{\PPi}{\mathcal{P}_i}
\newcommand{\EEi}{\mathcal{E}_i}
\newcommand{\beps}{\mbox{\boldmath $\epsilon$}}
\newcommand{\kmax}{k_{\mbox{\scriptsize max}}}

\section{Introduction}

The request for accurate and reliable numerical simulations of complex phenomena is
growing exponentially across a wide range of scientific and engineering applications.
In the automotive, aerospace or in the oil\&gas industry, the use of very large
computational grids is a common practice to faithfully reproduce intricate geometries
and guarantee a high precision level. The size of the problems
can easily reach several hundreds or even thousands of millions of unknowns up to the point that
High Performance Computing (HPC) becomes a necessity.
Despite of the differences in the problem at hand, there is always the need to discretize
the underlying Partial Differential Equations (PDEs) to approximate the continuous problem
in an algebraic system of equations whose solution is obtained numerically. In large scale
simulations, the solution of linear systems of equations is often, by far, the most
expensive part in terms of computational resources of the entire simulation process,
taking ofen more than 90\% of the total time.

The design and implementation of numerical software for the solution to large scale
linear systems on massively parallel computers is not straightforward and often
requires several person-years of development and testing to get a reliable code with
satisfactory performance. For many years now, US national labs have been provided
the scientific community with open source software for the solution of linear (and
non-linear) systems of equations on distributed memory systems~\cite{PETSc19,Hypre02,Trilinos20}.
Parallelism in these packages is exploited mainly through MPI implementing suitable
Krylov subspace methods preconditioned by Algebraic Multigrid~\cite{XuZik17} (AMG) techniques.
It has been shown in a number of pubblications that, provided a large enough problem,
the resulting codes may scale up to several thousands of CPU cores.

However, in recent years scientific computing is undergoing a major evolution.
Due to the end of the golden age ruled by Moore's law, traditional multi-core CPU are, by now, supported  
by many-core accelerators featuring thousands of simple computing cores. This change is clearly reflected in the Top 500 ranking
\cite{top500} where the first positions are occupied by systems equipped
with many-core accelerators. In particular, Graphics Processing Units, devices originally proposed for their efficiency at manipulating
computer graphics and image processing, are increasingly employed in general purpose numerical computing applications.

Adapting well-established software solutions, originally designed for traditional CPUs, to the use of accelerators
may be not an easy task up to the point of requiring a complete rethinking of the whole implementation.
First attempts to port linear solvers to the GPU hardware were based on a simple straightforward
approach \cite{LiSaa13b,LiSaa13a}. More recently, also some AMG based
solvers entirely running on a single or multiple GPUs have been proposed exhibiting a
very good performance on classic linear algebra problems
\cite{BelDalOls12,GanEslZha14,NauArsCasCohDem15,BerDamPas19}.
However, in challenging real world problems such as those arising from structural
mechanics or fluid flow in highly heterogeneous formations, standard AMG solvers
may be slow to converge or even fail, so that more advanced approaches are needed.
In particular, the use of powerful smoothers based on approximate inverses
can be of great help as shown, for instance, in \cite{PalFraJan19,FraMagMazSpiJan19}.

Approximate inverse preconditioners have been thoroughly investigated in the end of the last
century by several authors \cite{BenMeyTum96,GroHuc97,ChoSaa98}. Their main advantage
with respect to other preconditioners is their practically ideal suitability for being applied in parallel
as they involve only sparse matrix-by-vector products. Parallel processing is
more difficult to exploit in the set-up stage, but there are several approximate inverse
variants able to perform well even on parallel computers. See, e.g.~\cite{BenTum99,Ben02}
for a thorough review. Approximate inverse preconditioners are becoming popular again
due to their potential in concurrent computations especially on GPU hardware
\cite{XuDinFanChe11,BerFil14,AnzHucBraDon18,GuiRenJia19,MouGraFil20a,MouGraFil20b}.

In this work, we focus on the Factored Sparse Approximate Inverse (FSAI) preconditioner
originally introduced in \cite{KolYer93} and further developed over the years by other
authors \cite{JanFerSarGam15}. FSAI is particularly attractive due to its natural
parallelism in both set-up and application and has already been successfully ported
on GPU in both its static \cite{BerBisFanJan16} and dynamic \cite{BerCarFraJan19}
pattern variants using CUDA, the parallel computing platform and programming model developed by NVIDIA for general computing on graphical processing units.
The main novelty presented here is a highly optimized hybrid implementation (CUDA+MPI) that takes advantage of multiple GPU cards
to speed-up the solution of very large linear systems arising from the simulation of a number of real-world phenomena.

Due to recent changes in CUDA we had to modify some of the kernels
developed in the past and to tune the new distributed implementation
for exploiting at its best the state-of-the-art interconnection
technology we use.

The paper is organized as follows. In Section~\ref{theory}, we provide some basic theory
about the aFSAI preconditioner that we present. In Section~\ref{sec:GPU_FSAI}, we briefly
recall the main GPU kernels used for aFSAI set-up and their update to the most recent
CUDA version, whereas in section~\ref{DistrMem} we describe the issues related to the
distributed memory implementation. Extensive numerical experiments are presented in
section~\ref{NumRes}. Finally, we close the paper drawing some concluding remarks.

\section{Factored Sparse Approximate Inverse with dynamic pattern selection}
\label{theory}

\begin{figure}
%\centerline{\includegraphics[height=8cm]{./Figures/Pattern_EPS}}
\centerline{\psfig{figure=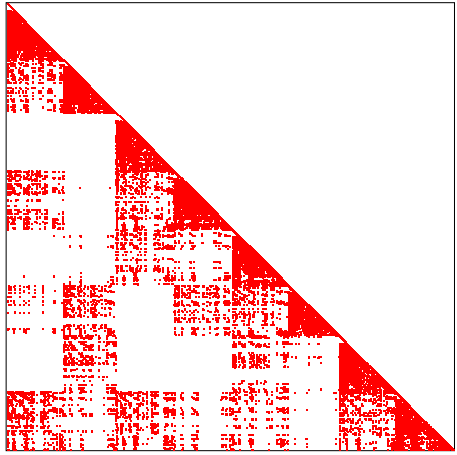,width=6cm}}
\caption{Example of a typical lower triangular non-zero pattern $\mathcal{S}$.}
\label{fig:patt}
\end{figure}

The FSAI preconditioner has been originally introduced for SPD matrices
in~\cite{KolYer93}, with the aim of directly approximating the inverse of
$A$ as the product of two triangular factors:
\begin{equation}
M^{-1} = G^T G \simeq A^{-1}
\label{prec}
\end{equation}
In the equation \ref{prec}, $G$ is a lower triangular matrix whose entries are computed
by minimizing the following Frobenius norm:
\begin{equation}
\| I - G L \|_F
\label{Froben_1}
\end{equation}
over the set $\mathcal{W}_{\mathcal{S}}$ of matrices having a prescribed lower triangular
non-zero pattern $\mathcal{S}$, as the one depicted in Figure~\ref{fig:patt}.
The matrix $L$, explicitly appearing in (\ref{Froben_1}), is the exact lower triangular
factor of $A$ and it is not actually needed in the computation of FSAI, since it
disappears during the minimization process as shown for instance in \cite{KolYer93}.
The unknown $G$ entries, ${[G]}_{ij}$, are computed by solving the componentwise system:
\begin{equation}
[G A]_{ij} = \left\{
\begin{array}{cccc}
0 & & i \neq j, & (i,j) \in {\mathcal S} \\
{[L]}_{ii} & & i = j  \\
\end{array}
\right.
\label{syst_2}
\end{equation}
obtained through differention of (\ref{Froben_1}) with respect to ${[G]}_{ij}$ and setting
it equal to zero. In the equation above, the symbol $[\cdot]_{ij}$ in (\ref{syst_2})
is used to indicate the entry in row $i$ and column $j$ of the matrix between square
brackets. Since ${[L]}_{ii}$, the $i$-th diagonal element of $L$, is unknown, we replace
it in eq. (\ref{syst_2}) by 1. As a consequence, in place of $G$, we compute the matrix
$\widetilde{G}$ by solving:
\begin{equation}
[\widetilde{G} A]_{ij} = \delta_{ij}
\label{syst_3}
\end{equation}
with $\delta_{ij}$ the Kronecker delta.

From a practical viewpoint, setting up the $i$-th row of $\widetilde{G}$,
say $\widetilde{\bg}_i^T$, requires:
\begin{itemize}
   \item first, the definition of the set $\PPi$ collecting all the column indices that
   belong to the $i$-th row of $\mathcal{S}$, that is:
   \begin{equation}
   \PPi = \{ j : (i,j) \in \mathcal{S} \}
   \end{equation}
   \item then, the gathering of the dense matrix $A[\PPi,\PPi]$ formed with the entries
   of $A$ having row/column indices in $\PPi$;
   \item finally, solving the linear system:
   \begin{equation}
   A[\PPi,\PPi] \widetilde{G}[i,\PPi]^T = \be_{m_i}
   \label{DensSys}
   \end{equation}
   with $\widetilde{G}[i,\PPi]$ being the dense vector containing the non-zero entries of
   $\widetilde{\bg}_i$ and the right-hand side $\be_{m_i}$ given by the last vector
   of the canonical basis of $\mathbb{R}^{m_i}$, with $m_i = |\PPi|$.
\end{itemize}

\begin{figure}
%\centerline{\includegraphics[scale=0.7]{./Figures/Schema_Red}}
\centerline{\psfig{figure=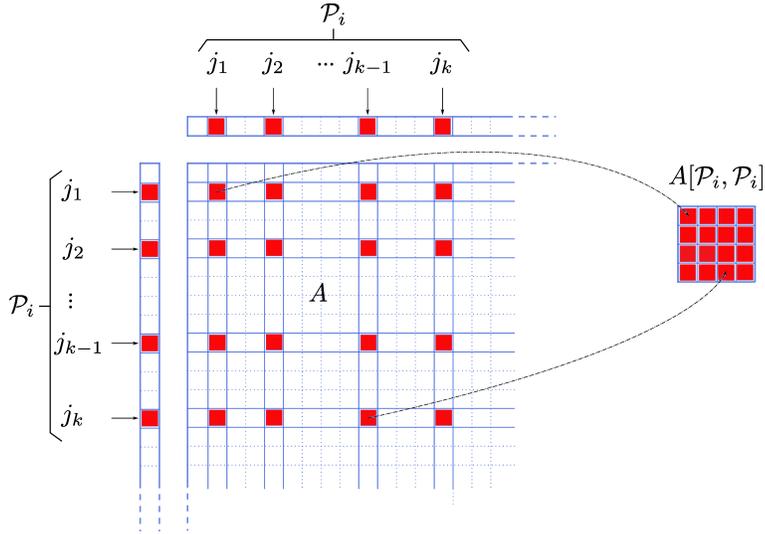,width=10cm}}
\caption{Schematic representation of the dense linear system gathering for a given set $\PPi$
(with cardinality $k$).}
\label{fig:gather}
\end{figure}

Figure \ref{fig:gather} gives an idea of the gathering process used to collect the dense
linear systems $A[\PPi,\PPi]$ given a set $\PPi$.

A more practical, although mathematical equivalent, way to implement the FSAI computation,
see~\cite{JanFerSarGam15}, consists in assuming $\widetilde{G}$ unitary diagonal and solving
the linear system:
\begin{equation}
A[\overline{\PPi},\overline{\PPi}] \widetilde{G}[i,\PPi]^T = -A[\overline{\PPi},i]
\label{blocksys}
\end{equation}
where $\overline{\PPi} = \PPi \setminus i$ and $\widetilde{G}[i,\PPi]$ contains
the {\em off-diagonal} non-zero entries of $\widetilde{G}$.
A diagonal scaling is finally applied to $\widetilde{G}$ in order to guarantee that all
the diagonal entries of the preconditioned matrix are unitary:
\begin{equation}
\ddiag(D_G \tg A \tg^T D_G) = \ddiag(GAG^T) = I
\label{diagScal}
\end{equation}
where $\ddiag(\cdot)$ is the operator returning the diagonal matrix having the diagonal
of its argument as entries, $I$ is the identity matrix and the $D_G$ entries are
given by:
\begin{equation}
D_G[i,i] = \frac{1}{A[i,i] - \widetilde{G}[i,\PPi] A[\overline{\PPi},\overline{\PPi}]
\widetilde{G}[i,\PPi]^T}
\label{diagScal2}
\end{equation}
Condition (\ref{diagScal}) ensures that $G$, over all the matrices
$B \in \mathcal{W}_{\mathcal{S}}$, is the unique one minimizing the Kaporin number
of the preconditioned matrix:
\begin{equation}
\kappa = \frac{\displaystyle \frac{1}{n} \tr(G A G^T)}{\displaystyle \dt(G A G^T)^{\frac{1}{n}}}
\end{equation}
which gives a measure of the PCG convergence rate~\cite{Kap94}.

The FSAI preconditioner for an SPD matrix is breakdown-free and possesses an extremely
high degree of parallelism in both construction and application to a vector.
However, to be competitive with other preconditioners, FSAI needs to be sparse
with a non-zero pattern composed by significant entries of the true inverse factor
of $A$. The appropriate a priori choice of $\mathcal{S}$ is unfortunately very
difficult thus making the original FSAI unpractical in ill-conditional cases.

A better way to compute FSAI is by choosing $\mathcal{S}$ {\em adaptively} while
computing $\widetilde{G}$ \cite{JanFer11}. The basic concept in the adaptive computation
of FSAI (aFSAI from now on) is the improvement of the quality of a given initial
factor $G_0$, already satisfying equations~(\ref{blocksys}) and~(\ref{diagScal}),
by extending its pattern with those entries that mostly contribute in reducing the Kaporin
number of $G_0 A G_0^T$.
Let us define $\mathcal{S}_0$ the non-zero pattern of $G_0$ and assume the identity
matrix as our default initial guess.
Writing explicitly the Kaporin number $\kappa$ of $G_0 A G_0^T$, it gives:
\begin{equation}
\kappa = \frac{\displaystyle \frac{1}{n} \tr(G_0 A G_0^T)}
{\displaystyle \dt(G_0 A G_0^T)^{\frac{1}{n}} } =
\frac{\displaystyle \frac{1}{n} \tr(D_{G_0} \tg_0 A \tg_0^T D_{G_0})}
{\displaystyle \dt(D_{G_0} \tg_0 A \tg_0^T D_{G_0})^{\frac{1}{n}} }
\label{kap}
\end{equation}
and, recalling that $G_0 A G_0^T$ has unitary diagonal entries and
$\dt(\widetilde{G}_0) = 1$ by construction, it follows that:
\begin{equation}
\kappa = \frac{1}{\displaystyle \dt(A)^{\frac{1}{n}}} \;
\frac{1}{\displaystyle \dt(D_{G_0})^{\frac{2}{n}}}=
\frac{\displaystyle \dt\left[\ddiag(\widetilde{G}_0 A \widetilde{G}_0^T)\right]^{\frac{1}{n}}}
{\displaystyle \dt(A)^{\frac{1}{n}}}
\label{kap1}
\end{equation}
Denoting by $\tbg_{0,i}^T$ the $i$-th row of $\widetilde{G}_0$, we can write the
numerator of (\ref{kap1}) as:
\begin{equation}
\dt\left[\ddiag(\widetilde{G}_0 \widetilde{A} \widetilde{G}_0^T)\right] =
\prod_{i=1}^n \tbg_{0,i}^T A \tbg_{0,i} = \prod_{i=1}^n \psi_{0,i}
\label{kap_min}
\end{equation}
having defined $\psi_{0,i} = \tbg_{0,i}^T A \tbg_{0,i}$.
With the above definiton, we can write a simplified expression for
the Kaporin conditioning number of the preconditioned matrix:
\begin{equation}
\kappa = \left( \frac{ \prod_{i=1}^n \psi_{0,i} }{\dt(A)} \right)^{\frac{1}{n}}
\label{kap2}
\end{equation}
which can be differentiated to obtain the gradient of $\kappa$ with respect to
$\widetilde{G}_0$. Choosing the positions of the gradient corresponding to its largest
entries in absolute value, allows for finding an augmented pattern $\cal{S}$$_{1}$ giving
a large reduction of $\kappa$ in equation~(\ref{kap2}).
Since each $\psi_{0,i}$ in~(\ref{kap2}) does not depend on any other row of
$\widetilde{G}_0$ except the $i$-th one, the pattern expansion can be carried out on each
row concurrently.
The $\psi_{0,i}$ gradient, denoted as $\nabla \psi_{0,i}$, is obtained by collecting
its partial derivatives with respect to the components of $\tbg_{0,i}$:
\begin{equation}
\frac{\partial \psi_{0,i}}{\partial {[\widetilde{G}_0]}_{ij}} = 2
\left( \sum_{r = 1}^n a_{jr}\, {[\widetilde{G_0}]}_{ir} + a_{ji} \right),
\qquad \forall j = 1,\dots,i-1,
\label{grapsi}
\end{equation}
with the computational cost for its evaluation consisting in a single sparse-matrix
by sparse-vector product.
A new pattern $\overline{\PP}^1_i$ is obtained by enlarging $\overline{\PP}^0_i$ with the
$s$ positions corresponding to the largest entries of $\nabla \psi_{0,i}$ in absolute
value and finally $\widetilde{\bg}_{1,i}$ is computed by solving~(\ref{blocksys}),
after replacing the sub/superscript ``0'' with ``1''.
Repeating the above procedure $\kmax$ times allows to find all the rows of $\widetilde{G}_2$,
$\widetilde{G}_3$, $\ldots$, up to $\widetilde{G}_{\kmax}$. It is also possible to monitor
the preconditioner quality by computing the value of $\psi_{k,i}$.
Hence it is possible to stop the adaptive procedure independently for every row when
\begin{equation}
\frac{\psi_{k,i}}{\psi_{0,i}} = \frac{\tbg_{k,i} A \tbg_{k,i}^T}{\tbg_{0,i} A \tbg_{0,i}^T}
\le \varepsilon,
\label{exit_crit}
\end{equation}
being $\varepsilon$ a user-specified tolerance.
Once the approximation of $\tg$ is satisfactory, a diagonal scaling is applied,
to ensure a unitary diagonal for $GAG^T$.

Three user-specified parameters are generally used to control the aFSAI quality:
\begin{enumerate}
\item $\kmax$, the maximum number of steps of the iterative procedure for each row;
\item $s$, the number of new positions added to the non-zero pattern of
each row in a single step;
\item $\varepsilon$, the relative tolerance on the Kaporin number reduction used to stop
the procedure.
\end{enumerate}

\newcommand{\txr}[1]{\textcolor{red}{#1}}

\section{GPU kernels for aFSAI set-up and application}
\label{sec:GPU_FSAI}

A detailed description of the CUDA kernels in charge of the most demanding, from the computational viewpoint, phases of the aFSAI set-up
is available in our previous works \cite{BerBisFanJan16,BerCarFraJan19}. Here we just recall that there are three main kernels for
\begin{enumerate}
\item the computation of $\nabla \psi_{i}$, i.e. of the {\em Kaporin gradient}, and selection of its most relevant components;
\item the collection of $A$ entries to form $A[\overline{\PP}_i,\overline{\PP}_i]$ and
$A[\overline{\PP}_i,i]$, that we call {\em systems gather};
\item the solution to the sequence of collected dense SPD systems, carried out by means of {\em batch Cholesky
decompositions}.
\end{enumerate}
As we showed in~\cite{BerCarFraJan19}, the computation of the Kaporin gradient of the input matrix is the
part of the preconditioner set-up that requires more time. Although the algorithm for the evaluation of the Kaporin gradient
remained the same, we had to modify its CUDA implementation due to the changes introduced in the {\em warp-level} primitives
starting on the CUDA toolkit version 9.0. Since CUDA has been around for a while, we are not going to describe it in detail.
However, we believe it is important to recall that, in the CUDA, GPUs execute warps of 32 parallel threads according to an execution model that NVIDIA calls 
Single Instruction Multiple Thread (SIMT), in other words, a variant of the SIMD (Single Instruction, Multiple Data) model, which is one of the four classes defined by
the classic Flynn's taxonomy. CUDA toolkits prior to version 9.0 relied on the implicit assumption that the threads within
a warp worked in a fully synchronous way. If this is not true, a program may show unexpected side effects up to the point of being unreliable.
The point can be illustrated in a simple case based on the \verb|__ballot(|{\em predicate}\verb|)| primitive that returns an \verb|unsigned int|
whose $n^{th}$ bit is set if-and-only-if {\em predicate} evaluates as \verb|true| for the $n^{th}$ thread of the warp. In the
fragment of code shown in listing \ref{listing:ballot}, the CUDA compiler and the hardware should try to re-converge the threads
right after the \verb|if/else| block for better performance. But this re-convergence is not guaranteed in the most recent versions of the CUDA toolkit.
Therefore, the \verb|ballot_result| variable may not contain the ballot result from {\em all} 32 threads.

\begin{lstlisting}[caption={Unsafe CUDA programming based on the implicit assumption that warp's threads run synchronously.},label=listing:ballot]
int result;
if (thread_id % 2) {
    result = foo();
} else {
    result = bar();
}
unsigned ballot_result = __ballot(result);
\end{lstlisting}
Starting on version 9.0, up to version 10.0 of the CUDA toolkit the legacy warp-level primitives worked synchronously (albeit with a deprecation warning at compilation time).
But, starting on version 10.1, the only alternative to obtain the expected behaviour, is to employ an explicit control on the threads that participate in warp operations by using the new form of the warp-level
primitives.
The set of threads that participates in each primitive is specified by means of a 32-bit mask, which is always the first argument in the new syntax of the warp-level primitives.
So, for instance, the new form of the \verb|__ballot()| primitive is \verb|__ballot_sync(mask,| {\em predicate}\verb|)|. 
All the participating threads are synchronized before the execution if they are not already synchronized. In the simple case reported in listing \ref{listing:ballot},
it is enough to replace the last line with \small
\verb|unsigned ballot_result = __ballot_sync(0xFFFFFFFF, result);|. \normalsize
Unfortunately there is not a general rule to determine what is the right value of the mask argument.
The set of threads to be included in the mask is determined by the program logic, and may depend on branch conditions or other information available only at execution time.
That is exactly the situation we had with the kernel for the evaluation of the Kaporin gradient that relies heavily on warp-level primitives like 
{\em shuffle} to exchange values stored in the registers among threads belonging to the same warp (as described in \cite{BerCarFraJan19} this is the technique that makes 
very efficient the selection of min and max values within that kernel).

\begin{figure*}[h]
\centering
\includegraphics[angle=-90,width=0.8\textwidth]{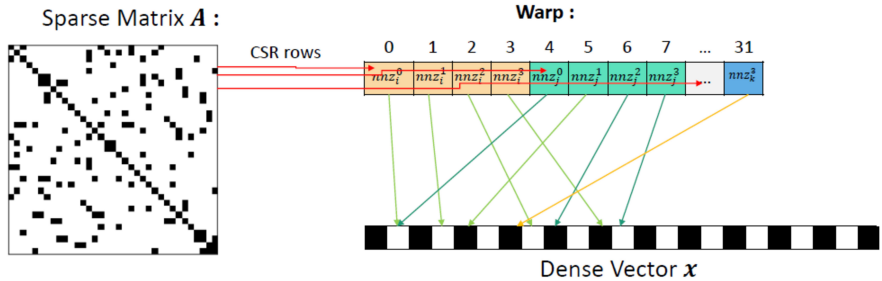}
\caption{Each {\em miniwarp} is in charge of a row of the matrix stored in CSR format.\label{fig-miniwarp}}
\end{figure*}

Besides the update of the kernel for the Kaporin gradient, for the
present work, which includes a new multi-GPU implementation of the
conjugate gradient algorithm, we developed also a new CUDA kernel for
the sparse-matrix-dense-vector product. To that purpose we employ a
group of threads for each row of the sparse matrix. Using a group instead of a single thread
makes it possible to exploit at its best the memory bandwidth of the GPU,
at least loading the elements of the matrix. We call the group of
threads in charge of each row \textit{miniwarp} in analogy with the
group of $32$ consecutive threads named \textit{warp} in the CUDA
jargon that we mentioned above. In particular, we
consider sets (having the same number of elements) of contiguous
threads with a possible cardinality of $2$, $4$, $8$, $16$ or
$32$. Those sets of threads are then concurrently mapped to different
data like in other warp-centric kernels that use a warp to manage a
block of contiguous data. The threads contained in a miniwarp are able
to perform cooperative computation by exploiting efficient and
fine-grained intra-warp communication primitives like the {\em shuffle}.
During the execution of the product, each row of the sparse matrix is
assigned to a single miniwarp; that is, multiple rows are concurrently
executed in the same full warp of $32$ threads. The size of a miniwarp
is dynamically dependent on the average number of nonzeroes per
row of the sparse matrix. The main advantage of this miniwarp variant is
that, for matrices with few nonzero entries per row, the number of
idle threads decreases. With a full warp (32 threads), if a row has, on average,
only $k < 32$ nonzero entries, there are, always on average, $32 -k$
threads that remain idle. Miniwarps reduce the difference
significantly by using a size that is much closer to the average
number of nonzero entries per row. Fig.~\ref{fig-miniwarp} shows how
the miniwarp works when applied to a matrix in CSR format in the case
of a sparse-matrix dense-vector multiplication.

The GPUs used during the development of the present work are
NVIDIA Tesla V100 (based on the Volta CUDA architecture) equipped with
80 Stream Multiprocessors for a total of 5120 cores and 16 GByte of
GDDR5 memory.

We show an example of the contribute of the three main kernels in the following
fragment of the CUDA profiler ({\tt nvprof}) output for a medium size matrix:

\begin{verbnobox}[\footnotesize]
   Time(%)      Time  Calls       Avg       Min       Max  Name
   59.82%   5.12072s    400  12.802ms  698.88us  26.824ms  kapGrad_merge
   27.27%   2.33450s     20  116.72ms  16.199ms  180.44ms  gatherFullSys
    4.96%   424.63ms     20  21.232ms  170.02us  76.376ms  cudacholdc
   ...
\end{verbnobox}

The three kernels take more than $90\%$ of the time. The remaining kernels are used
to check if the convergence criterion for the refinement has been reached and to build,
in the end, the preconditioner. None of them takes more than $0.5\%$ of the total
execution time and are all pretty simple.

\section{Distributed memory implementation}
\label{DistrMem}

In this section we give a detailed description of the distributed memory implementation
of the aFSAI preconditioner and its use in an iterative solver as the conjugate gradient.
To this aim, we firstly introduce the storage scheme adopted for the system matrix
because from this scheme descends the communication pattern in both
aFSAI set-up and Sparse-Matrix-by-Vector product (SpMV). We observe that the proposed 
storage scheme applies to the hybrid CPU-GPU implementation as well as to the pure CPU
implementation that we use as a baseline for our numerical experiments. In the
description that follows, we do not distinguish between the two implementations as the
only difference is in the few, though computationally intensive, specific kernels
introduced in section~\ref{sec:GPU_FSAI}.

For the communication among different computational nodes we use the Message Passing
Interface (MPI), whereas for intranode parallelization we use OpenMP for the CPUs and CUDA for the GPUs.

\subsection{Distributed Sparse (DSMat) and Dense (DDMat) Matrix Storage Schemes}

\begin{figure*}[h]
\centering
\includegraphics[angle=-90,width=0.4\textwidth]{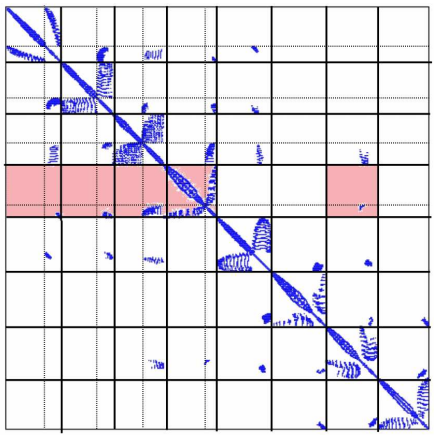}
\caption{Schematic representation of the DSMat matrix storage scheme using 8 MPI processes.
         The red colored blocks are assigned to process 3.\label{fig-dsmat}}
\end{figure*}

The storage scheme adopted to distribute the system matrix, referred to as DSMat, essentially
consists in subdividing the sparse matrix into $n_p$ horizontal stripes of consecutive rows,
where $n_p$ is the number of MPI processes. Each stripe is then further split vertically
into rectangular blocks by applying the same subdivision, as shown in Fig.~\ref{fig-dsmat}.
Finally, each matrix subblock is stored according to the popular Compressed Sparse Row (CSR) matrix
format.

To save memory and increase the overall efficiency, a local numbering is used within each CSR
subblock, that is rows and columns of block $IJ$ are numbered, in C notation, from 0 to $n_{I-1}$
and from 0 to $n_{J-1}$, respectively, with $n_I \times n_J$ the dimension of block $IJ$. This
expedient allows for the use of the more favourable 4-byte integers even for very large matrices
having billions of rows.
Each process stores the diagonal block, which is generally the most populated one, representing
the part of the matrix operating on the local unknowns, a list of ``Left'' and ``Right'' blocks
representing the communication with neighboring processes having lower or higher rank,
respectively. For instance, with reference to Fig.~\ref{fig-dsmat},
processor 3 stores the 5 blocks highlighted in red, of which 0, 1 and 2 are denoted as left
neighbors, 3 as diagonal block and 6 as right neighbor.

This blocked scheme used to store matrices is very effective in both the aFSAI computation and
SpMV product because, as described in details below, it allows the
use of non-blocking send/receive MPI primitives with a large overlap between computation and
communication.

Similarly to DSMat, we define the DDMat storage scheme for distributed dense matrices, where the
matrix is partitioned into $n_p$ horizontal stripes of consecutive rows. The entries of each
stripe are stored row-wise in memory to guarantee a better access in memory during multiplication
operations. This storage scheme is particularly advantageous for linear systems with multiple
right-hand-sides and eigenproblems. For standard linear systems with only one rhs vector, the
DDMat matrix simply reduces to one distributed vector.

\subsection{Distributed memory Sparse-Matrix-by-Vector Product (SpMV)}

\hfill \break
The SpMV product $A\by=\bz$ between a DSMat $A$ and a DDMat
$\by$ is the most expensive operation in any preconditioned iterative solver and on a
distributed memory computer consists of a $communication$ stage, mainly handled by the CPUs, and a
$computation$ stage, that can be performed either by the CPU or the GPU accelerator (or, in some implementations by both).

\begin{figure*}[h]
\centering
\includegraphics[angle=-90,width=0.4\textwidth]{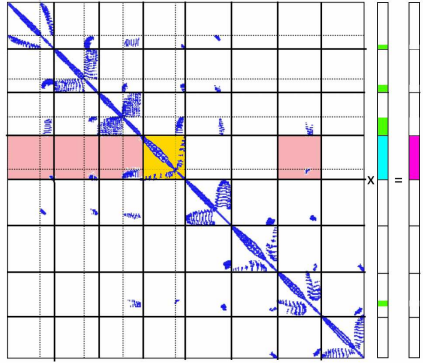}
\caption{Schematic representation of the SpMV product.\label{fig-spmv}}
\end{figure*}

With reference to Fig.~\ref{fig-spmv}, processor 3 gathers the $\by$ terms highlighted
in green and computes its stripe of $\bz$, highlighted in violet, as the sum of products
between the extra-diagonal CSR blocks, in red, and the entries received from neighbouring processes,
plus the contribution of the product between diagonal CSR blocks, in orange, and the stored
$\by$ terms, in light blue.
As detailed below, communication and computation can be conveniently overlapped thanks to
non-blocking communications and the DSMat storage scheme:
\begin{enumerate}
\item First, local elements of $\by$ to be sent are moved into a device buffer and copied
      on the host;
\item Then the host starts non-blocking send/receive communications to/from neighbouring processes;
\item in the meantime, the local component of $\bz$ is initialized with the product between
      the diagonal CSR block and the local part of $\by$;
\item each process sequentially tests all the incoming communications and as soon as some components
      of $\by$ are received, the corresponding buffer is copied onto the device and the local
      part of $\bz$ is updated with the product by the off-diagonal CSR block;
\end{enumerate}
It is easy to observe that the 3rd operation in the above list can be overlapped with the first two,
and also the product by the off-diagonal blocks partly hides the communication latency for the
transfer of non-local $\by$ components.

For an efficient execution, the SpMV algorithm requires a preliminary stage that we call {\em prepare},
during which each process communicates to its neighbors the list of $\by$ entries it needs to receive
to perform the product by the off-diagonal blocks. Since this information does not change during the 
iterative solution of the system, it needs to be exchanged just once before the calculation starts.

\subsection{aFSAI set-up on a distributed memory computer}

\begin{figure*}[h]
\centering
\includegraphics[angle=-90,width=0.4\textwidth]{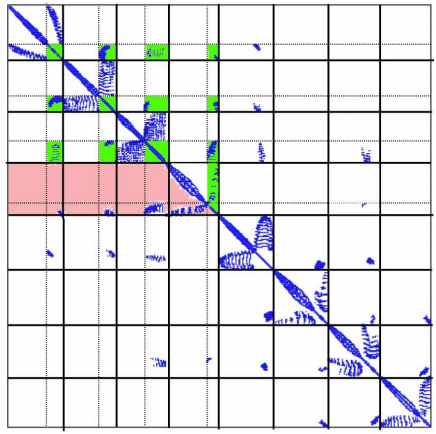}
\caption{Schematic representation of the block collection for process 3.
         Green colored blocks in the lower part of the matrix $A$ are those that need
         to be gathered for the assembly of $A_3$.
\label{fig-afsai_blocks}}
\end{figure*}

Before calling the specific GPU kernels (briefly described in Section~\ref{sec:GPU_FSAI}) for the aFSAI set-up on
every single device, each MPI task needs to collect some information from the others.
As the pattern of $G$ is computed dynamically during the set-up, in principle each task may need
matrix $A$ entries belonging to any other task with lower rank. Obviously, such a condition
cannot be satisfied in very large problems because the whole $A$ matrix cannot be stored
on a single node. On the other hand, preliminary tests showed that using only the information
from neighbouring processes, as in the communication pattern of SpMV, is not enough to build
a satisfactory preconditioner. So, we adopted an intermediate strategy. Let us denote by
$\widehat{A}$ the communication matrix of $A$, that is a boolean matrix of size
$n_p \times n_p$ with entries $\widehat{A}_{ij} \ne 0$ if and only if the corresponding block
$ij$ of $A$ contains at least one non-zero element. Then we impose that $\widehat{G}$, the
communication graph of $G$, is no larger than the lower triangular part of $\widehat{A}^k$ for
a small power $k$, typically in the range $1 \div 3$, which is simply computed through
a short sequence of small symbolic Sparse-Matrix-by-Sparse-Matrix products.
Once $\widehat{G}$ is known, we use it to guide the collection of information from other
MPI tasks.
With reference to Fig.~\ref{fig-afsai_blocks}, for instance, process 3 assembles locally
the matrix $A_3$, that it uses locally to compute its own stripe $G_3$ of $G$ by gathering
all the green blocks from its left-neighboring processes on $\widehat{G}$.
Thanks to the symmetry, only the communication of the green blocks in the lower part is needed
and the matrix in complete on site after a {\em transpose} operation.

Note that, again as a positive effect of the symmetry, each process knows in advance from the pattern of $\widehat{A}^k$
which blocks it has to receive and send to the other processes as well as the list of
processes it needs to communicate with. The gathering/scattering process consists of 2 stages
of communications, that can be partly overlapped by using non-blocking communication:
\begin{enumerate}
\item Each process sends to its right-neighboring processes (and receive from the left ones)
      the sizes of the blocks that have to be exchanged;
\item Each process sends to its right-neighboring processes (and receive from the left ones)
      the above set of blocks;
\end{enumerate}
We briefly note that, if a similar approach is extended to non-symmetric systems, each process
needs to communicate to the others also the list of required blocks, which is not known, a priori, as in
the symmetric case.

The communicated list of blocks can be conveniently received into a 2-dimensional array
data structure with size $(n_L+1)^2$, where $n_L$ is the number of left-neighboring
processes in $\widehat{G}$. Only the block-lower part of this structure is actually
filled with $({n_L}^2 + n_L)/2$ blocks received from the other processes plus the $(n_L+1)$
blocks already in place. Once all the communications have been completed, the 2-dimensional
array of CSR matrices is copied into the full (lower plus upper) CSR matrix $A_{I_p}$ that
is the main input for the aFSAI GPU kernel.
Note that, during this copy, these two large structures must coexist, becoming the main part of
the memory footprint of each process, that can be quantified as $\sim 1.6$ the size of $A_{I_p}$.

As already noted in \cite{BerCarFraJan19}, the aFSAI set-up can be safely performed in
single precision and, in the distributed memory context, this option is particularly
advantageous. In fact, beyond reducing processing time due to the use of single-precision
arithmetic, it also allows a smaller amount of data movement during the block-collection
stage and the host-device copies. Two casting operations are needed in this case.
The first one is performed on the matrix $A$ at the host level: $A$ is copied into a
new matrix $A_s$ using single precision floats for its entries, $A_s=\mbox{single}(A)$.
Then, the matrix $A_s$ is used for the aFSAI set-up, as described above.
Once $G_s$, $G$ with single precision entries, is computed, a second cast is performed at
the device level transforming single precision floats into doubles, $G=\mbox{double}(G_s)$.
This last casting operation can be conveniently overlapped with the host-device copy of
the double precision $A$, since only $A_s$ has been transferred before.

When both the double precision representation of $A$ and $G$ are present on the device,
the PCG algorithm is ready to start.

\section{Numerical Results}
\label{NumRes}

The aFSAI preconditioner is implemented as part of Chronos, a sparse linear algebra library
for preconditioners and Krylov subspace iterative methods specifically designed
for High-Performance computing systems~\cite{Chronos}.
For the evaluation of aFSAI performance, the Marconi100 supercomputer is used.
Marconi100, the most recent cluster installed at the Italian consortium for computing,
CINECA, is classified within the first ten positions of the TOP500 ranking. It is composed by
980 nodes based on the IBM Power9 architecture each equipped with two 16-cores
IBM POWER9 AC922 at 3.1 GHz processors and four NVIDIA Volta V100 GPUs with Nvlink 2.0 and
16GB of memory.
The benchmark matrices used for the numerical experiments are shown in Table~\ref{tab-matrices},
which provides the number of rows, $\mbox{n}(A)$, the number of non-zeroes, $\mbox{nnz}(A)$
and the field of application the matrices arise from. Most of the matrices have been collected
over years from internal research or cooperation with other authors
\cite{KorLuGul14,JhaRub11,Nat17,Nat20}.
In the remainder of this section, we first compare the efficiency of Chronos to that of the
well-known open source package PETSc~\cite{PETSc19} in CPU-only runs to draw a baseline for
the GPU experiments. Then, we compare GPU-accelerated aFSAI with the more traditional and
perfectly scalable Jacobi preconditioner on our benchmark set. Finally, we provide some
strong and weak scalability tests to prove the efficiency of the proposed implementation
while using large scale computational resources.

% Info sulle matrici utilizzate
\begin{table}
\begin{center}
\centering
{\small 
\begin{tabular}{lrrr}
\hline
Matrix Name & $\mbox{n}(A)$ & $\mbox{nnz}(A)$ & Application field \\
\hline
spe10     &   3,410,693  &     90,568,237 & 3D diffusion in heterogeneous media\\
geo4m     &   4,224,870  &    335,738,340 & 3D geomechanics \\
Finger    &   4,718,592  &     23,591,424 & 2D multiphase flow in porous media \\
guenda11m &  11,452,398  &    512,484,300 & 3D geomechanics \\
M10       &  11,593,008  &    940,598,090 & 3D mechanical \\
agg14m    &  14,106,408  &    633,142,730 & 3D mesoscale \\
M20       &  20,056,050  &  1,634,926,088 & 3D mechanical \\
geo61m    &  61,813,395  &  4,966,380,225 & 3D geomechanics \\
Pflow73m  &  73,623,733  &  2,201,828,891 & 3D diffusion in heterogeneous media\\
c4zz134m  & 134,395,551  & 10,806,265,323 & 3D Biomedicine \\
pois198m  & 198,076,032  &  1,384,390,392 & 3D diffusion in omogeneous media \\
\hline
\end{tabular} }
\end{center}
\caption{Benchmark matrices used in the numerical experiments. For each matrix, the size,
$\mbox{n}(A)$, the number of non-zeros, $\mbox{nnz}(A)$, and the application field are
provided.}
\label{tab-matrices}
\end{table}

% Benchmark spe-3m
\begin{figure}[htb!]
\begin{center}
\psfig{figure=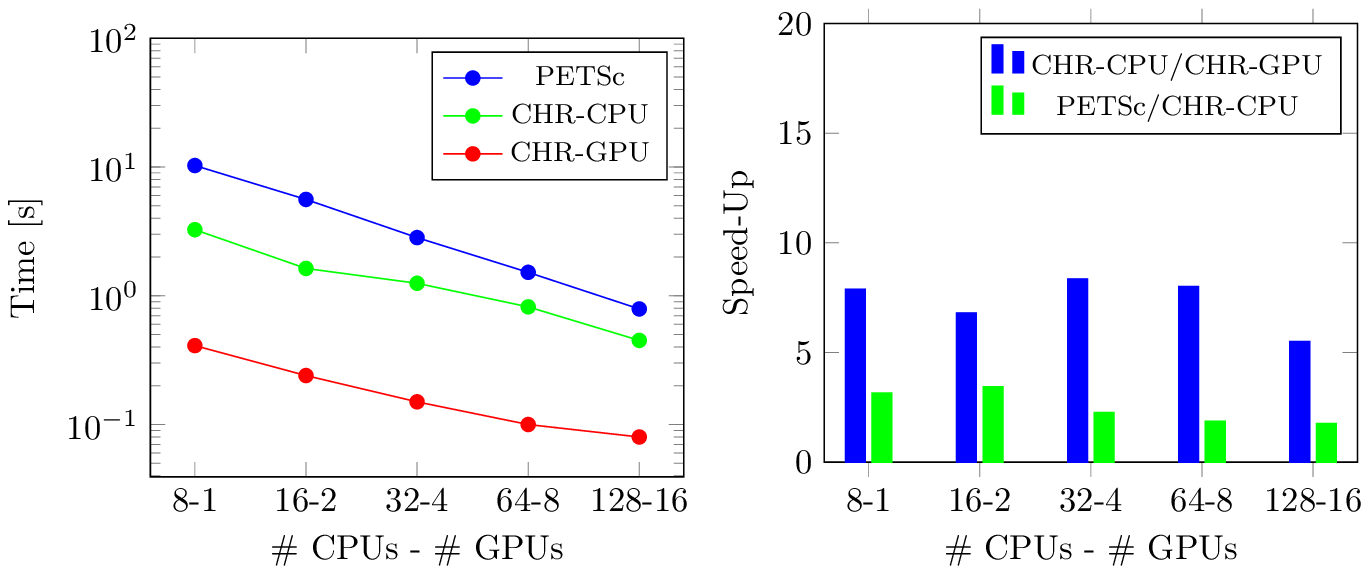,height=5cm}
\caption{Scalability and Speed-Up of Jacobi preconditioned CG on the
{\tt spe10} matrix. Left: total wall time in seconds for the execution of 100
iterations of Jacobi preconditioned CG. Right: Speed-Up of GPU-accelerated Chronos
(CHR-GPU) over pure CPU Chronos (CHR-CPU) (blue columns) and Speed-Up of Chronos over
PETSc in pure CPU runs (green columns). The computing resources vary from a quarter
node to 4 nodes.}
\label{fig-bench-spe}
\end{center}
\end{figure}

% Benchmark geo-4m
\begin{figure}[htb!]
\begin{center}
\psfig{figure=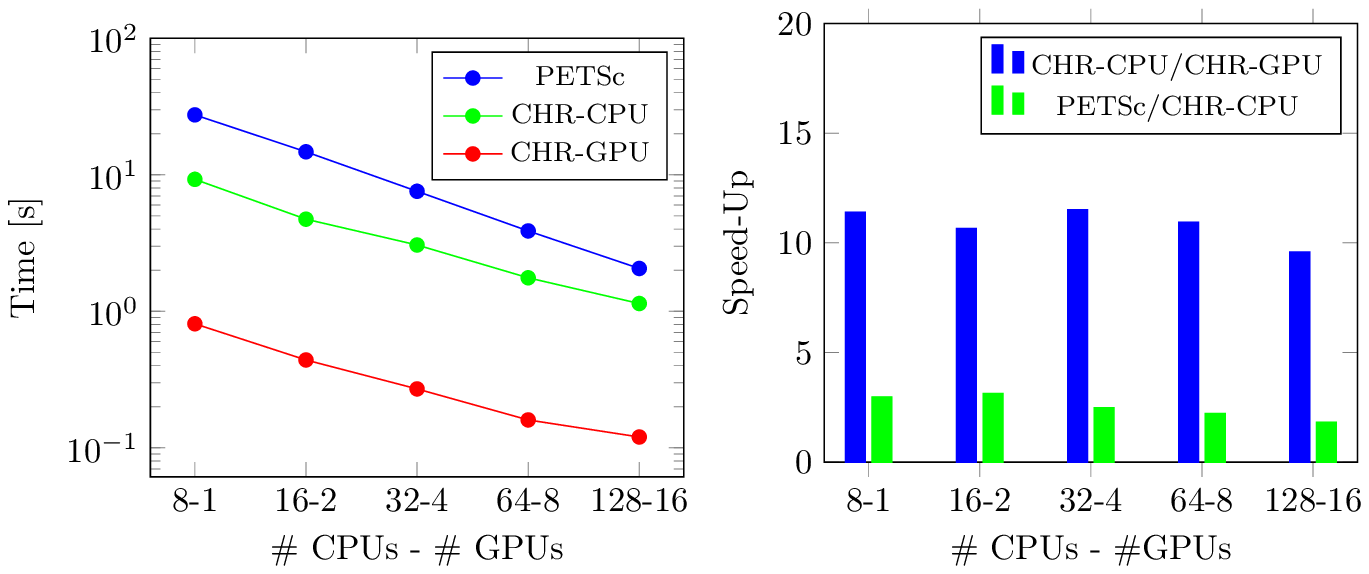,height=5cm}
\caption{Scalability and Speed-Up of Jacobi preconditioned CG on the
{\tt geo4m} matrix. Left: total wall time in seconds for the execution of 100
iterations of Jacobi preconditioned CG. Right: Speed-Up of GPU-accelerated Chronos
(CHR-GPU) over pure CPU Chronos (CHR-CPU) (blue columns) and Speed-Up of Chronos over
PETSc in pure CPU runs (green columns). The computing resources vary from a quarter
node to 4 nodes.}
\label{fig-bench-geo}
\end{center}
\end{figure}

% Benchmark agg-14m
\begin{figure}[htb!]
\begin{center}
\psfig{figure=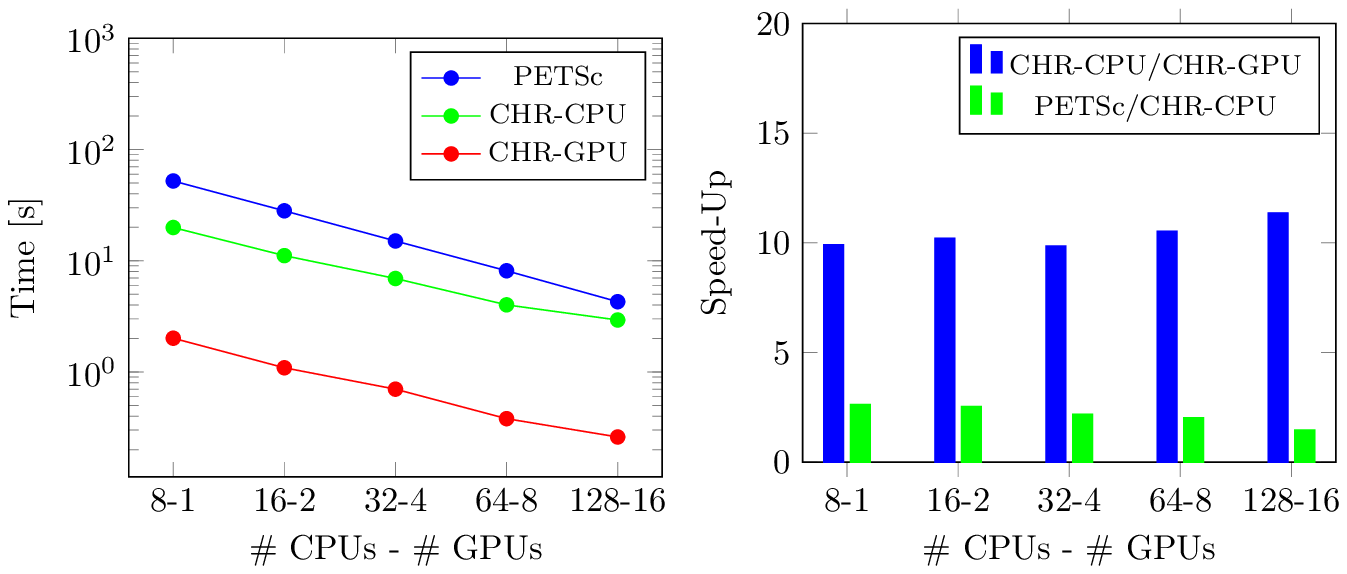,height=5cm}
\caption{Scalability and Speed-Up of Jacobi preconditioned CG on the
{\tt agg14m} matrix. Left: total wall time in seconds for the execution of 100
iterations of Jacobi preconditioned CG. Right: Speed-Up of GPU-accelerated Chronos
(CHR-GPU) over pure CPU Chronos (CHR-CPU) (blue columns) and Speed-Up of Chronos over
PETSc in pure CPU runs (green columns). The computing resources vary from a quarter
node to 4 nodes.}
\label{fig-bench-agg}
\end{center}
\end{figure}

\subsection{Effectiveness of the DSMat storage scheme in the sparse matrix by vector
product}

In this subsection, we show how the proposed DSMat storage scheme improves
the efficiency of the SpMV product in both CPU-only and GPU-accelerated runs.

To this aim, we first compare the time required by Chronos and PETSc to
perform 100 iterations of Conjugate Gradient (CG) preconditioned with Jacobi using only CPUs. The choice of Jacobi
is dictated by the fact that, since its implementation is straightforward, these tests 
allow us to directly evaluated the SpMV efficiency.
Figures \ref{fig-bench-spe} to \ref{fig-bench-agg} show the comparison between PETSc,
in pure CPU mode, and Chronos, in both CPU-only and GPU-accelerated mode, by increasing
the computing resources from a quarter of node (8 CPUs and 1 GPU) to 4 nodes
(128 CPUs and 16 GPUs) on the test matrices {\tt spe10}, {\tt geo4m} and {\tt agg14m}.
The speed-up of Chronos CPU (CHR-CPU) over PETSc is on average 2.39 and reaches a
maximum value of about 3.15. As expected, smaller speed-up values (about 2.00) are
obtained as the computing resources increase due to the communication overhead becoming
significant on the overall run time.
The GPU acceleration leads to a speed-up of about 10 for {\tt agg14m} and {\tt geo4m}
matrices while it reaches a smaller value around 6 for {\tt spe10}.
This lower performance is easily explained by considering the lower number of non-zeroes
per row characterizing {\tt spe10}, which induces a worse operation over communication
ratio. The use of GPU-accelerated Chronos allows an overall speed-up over CPU-only PETSc
in the range from 15 to 25 on the performed test cases.

\subsection{Effectiveness of aFSAI over Jacobi}

\begin{figure}[htb!]
\begin{center}
\psfig{figure=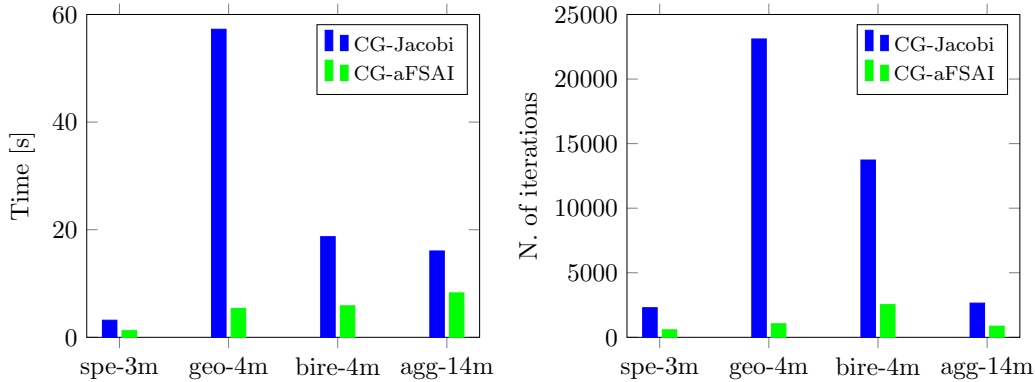,height=5cm}
\caption{Comparison between the GPU-accelerated CG preconditioned with Jacobi and
aFSAI using a single node of Marconi100.
Left: Total solution time. Right: Number of iterations to converge.}
\label{fig-jacVSfsai}
\end{center}
\end{figure}

The aFSAI preconditioner is much more demanding than Jacobi in term of both
implementation effort and set-up time. However, the adoption of aFSAI in real world
problems is fully justified by its superior effectiveness in accelerating CG convergence.
Figure~\ref{fig-jacVSfsai} shows the comparison between the total solution time, including
also the set-up time, and the number of iterations necessary to reduce the initial residual of
8 orders of magnitude using CG preconditioned with Jacobi and aFSAI.
We only consider the 4 smallest matrices of our set, as their solution cost though Jacobi
would have been prohibitive with the largest ones.
The experiments were run using the GPU-accelerated Chronos library on a single Marconi100
node and clearly show the ability of aFSAI in reducing the number of iterations. In all
the tests aFSAI outperfoms Jacobi by a factor of at least 2 in the worst case up to 10
in the most ill-conditioned {\tt geo4m}.

% Tempo di soluzione Jacobi VS aFSAI con 1 nodo di m100

\subsection{Scalability of the GPU-accelerated adaptive FSAI}

\begin{figure}[htb!]
\begin{center}
\psfig{figure=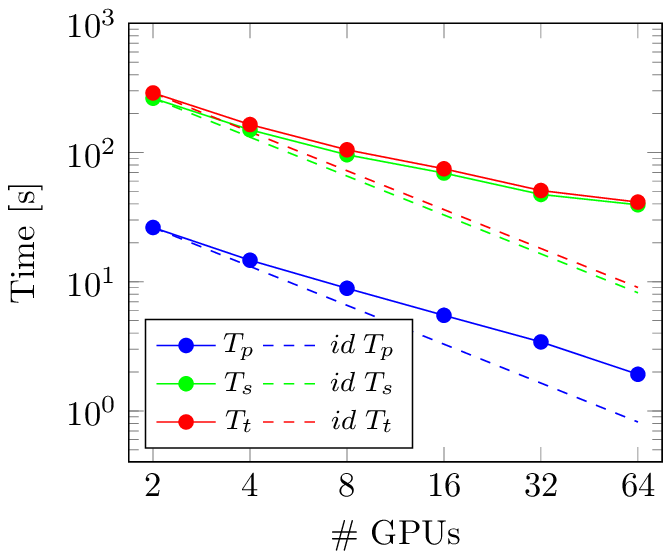,height=5cm}
\psfig{figure=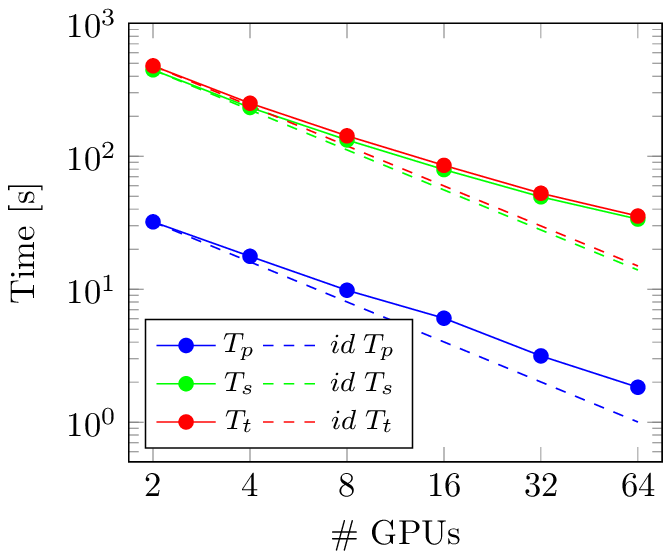,height=5cm}
\caption{Strong scalability of aFSAI CG on matrices {\tt guenda11m} (left) and {\tt M10} (right).
aFSAI set-up time $T_p$, CG iteration time $T_s$ and
overall time $T_t=T_p+T_s$ vs. number of GPUs.}
\label{fig-scal-guenda}
\end{center}
\end{figure}

\begin{figure}[htb!]
\begin{center}
\psfig{figure=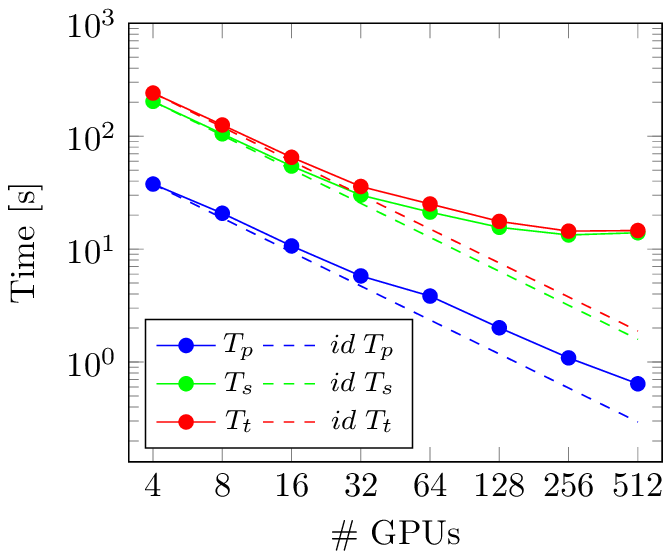,height=5cm}
\psfig{figure=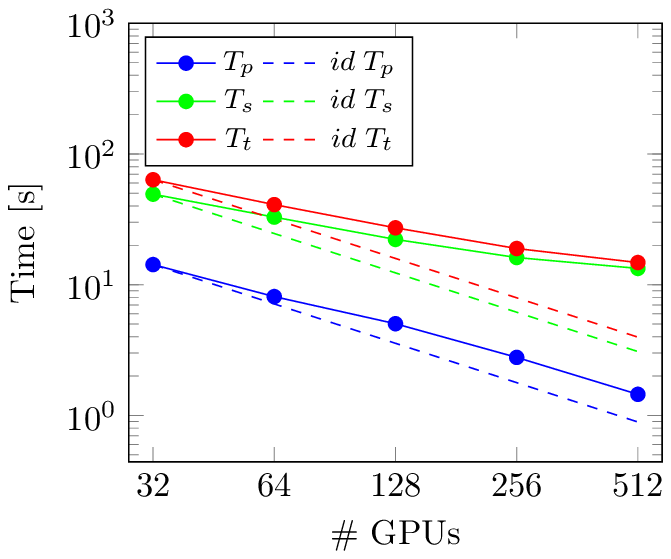,height=5cm}
\caption{Strong scalability of aFSAI CG on matrices {\tt M20} (left) and
{\tt geo61m} (right). aFSAI set-up time $T_p$, CG iteration time $T_s$ and
overall time $T_t=T_p+T_s$ vs. number of GPUs.}
\label{fig-scal-kor20}
\end{center}
\end{figure}

\begin{figure}[htb!]
\begin{center}
\psfig{figure=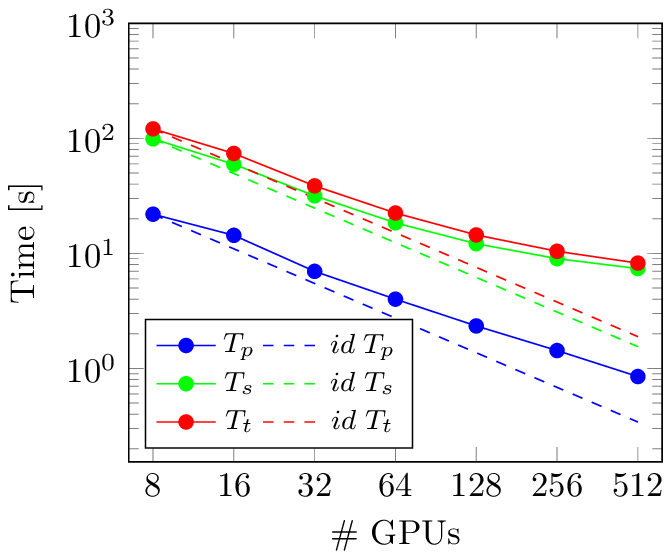,height=5cm}
\psfig{figure=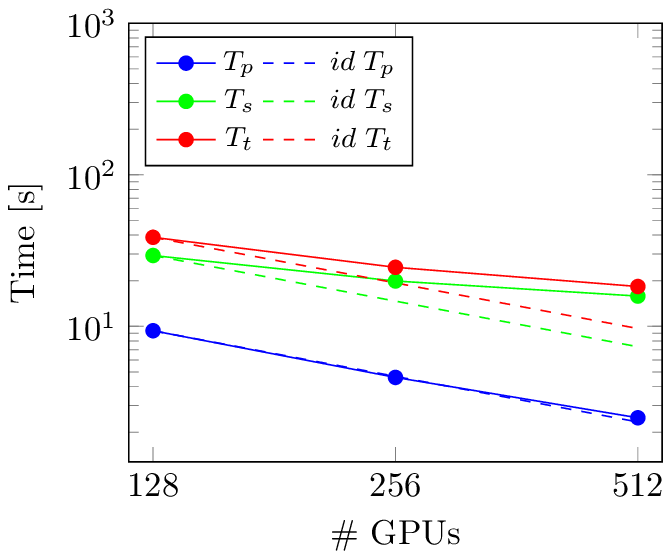,height=5cm}
\caption{Strong scalability of aFSAI CG on matrices {\tt Pflow73m} (left) and 
{\tt c4zz134m} (right).
aFSAI set-up time $T_p$, CG iteration time $T_s$ and
overall time $T_t=T_p+T_s$ vs. number of GPUs.}
\label{fig-scal-kutei}
\end{center}
\end{figure}

\begin{figure}[htb!]
\begin{center}
\psfig{figure=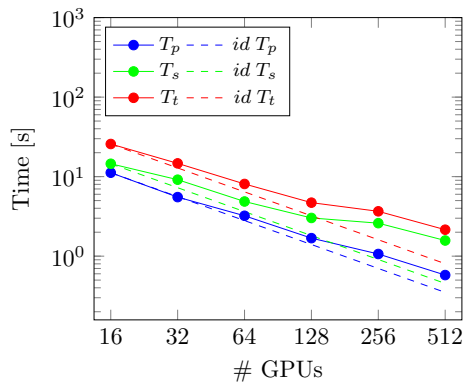,height=5cm}
\caption{Strong scalability of aFSAI CG on matrix {\tt pois198m}.
aFSAI set-up time $T_p$, CG iteration time $T_s$ and
overall time $T_t=T_p+T_s$ vs. number of GPUs.}
\label{fig-scal-poisson}
\end{center}
\end{figure}

\begin{figure}[htb!]
\begin{center}
\psfig{figure=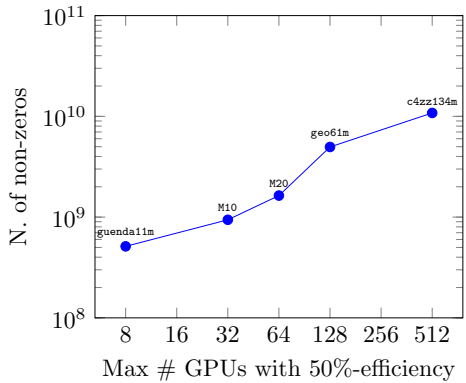,height=5cm}
\caption{Maximum number of GPUs that guarantee a parallel efficiency of, at least, 50\% for test matrices arising from structural mechanics.}
\label{fig-eff}
\end{center}
\end{figure}

In this last subsection we analyze the strong and weak scalability of aFSAI preconditioned
CG. We consider three different times to evaluate scalability: the preconditioner set-up
time, $T_p$, the iteration time, $T_s$, and total time, $T_t = T_p + T_s$.
We ran strong scalability tests on the 8 largest matrices of Table~\ref{tab-matrices},
by varying the number of GPUs from the minimum necessary to store the matrix and the
preconditioner to a maximum of 512, corresponding to 128 nodes.
Figures \ref{fig-scal-guenda} to \ref{fig-scal-poisson} provide execution times and
parallel efficiency vs the number of GPUs.
In all tests both the set-up and solution time decrease inversely proportional
to the computing resources with an almost ideal behaviour.
The set-up stage is affected less than the iteration stage by efficiency reduction
as communications have a lower impact in the preconditioner set-up.
Regarding the parallel efficiency, Figure \ref{fig-eff} shows the maximum number of GPUs such that the efficiency is, at least, 50\% for the matrices arising from structural
mechanics. We focused on this particular subset of matrices because all of them consists
of more than 10 millions rows and have a similar number of non-zeroes per row ranging
from 45 to 80.
As expected, the computational resources that can be used at the same efficiency
increase with the number of non-zeroes in the matrix.

% Scalabilità guenda-11m

% Scalabilità koric-20m

% Scalabilità kutei-73m

% Scalabilità poisson-198m

% Efficiency

% Weak scalability
\begin{figure}[htb!]
\begin{center}
\psfig{figure=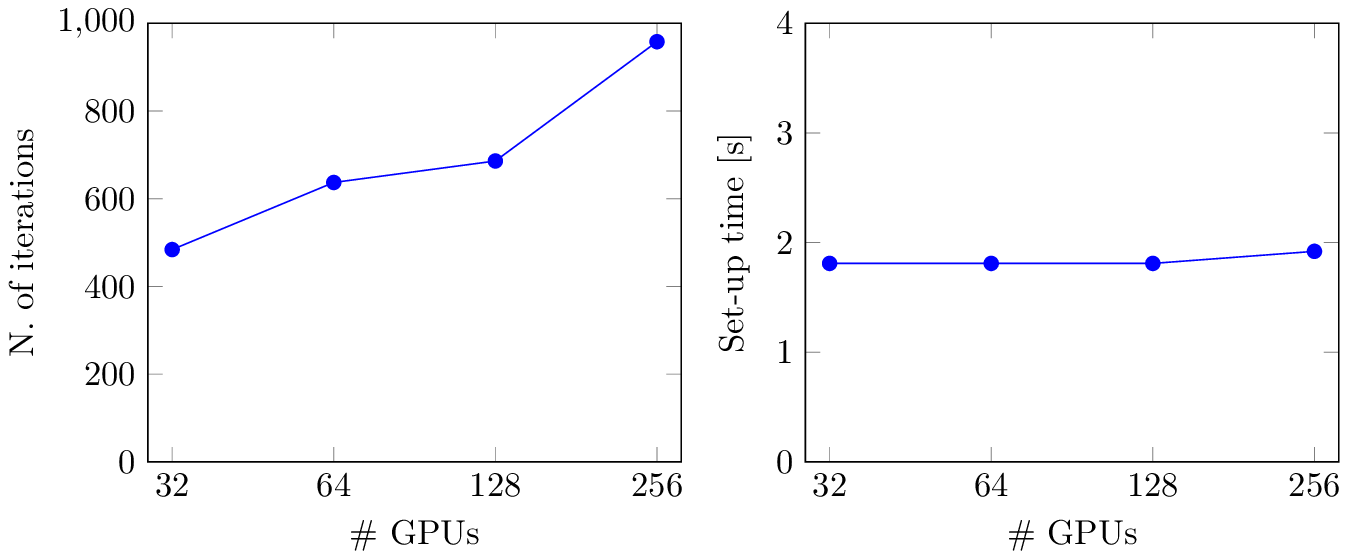,height=5cm}
\psfig{figure=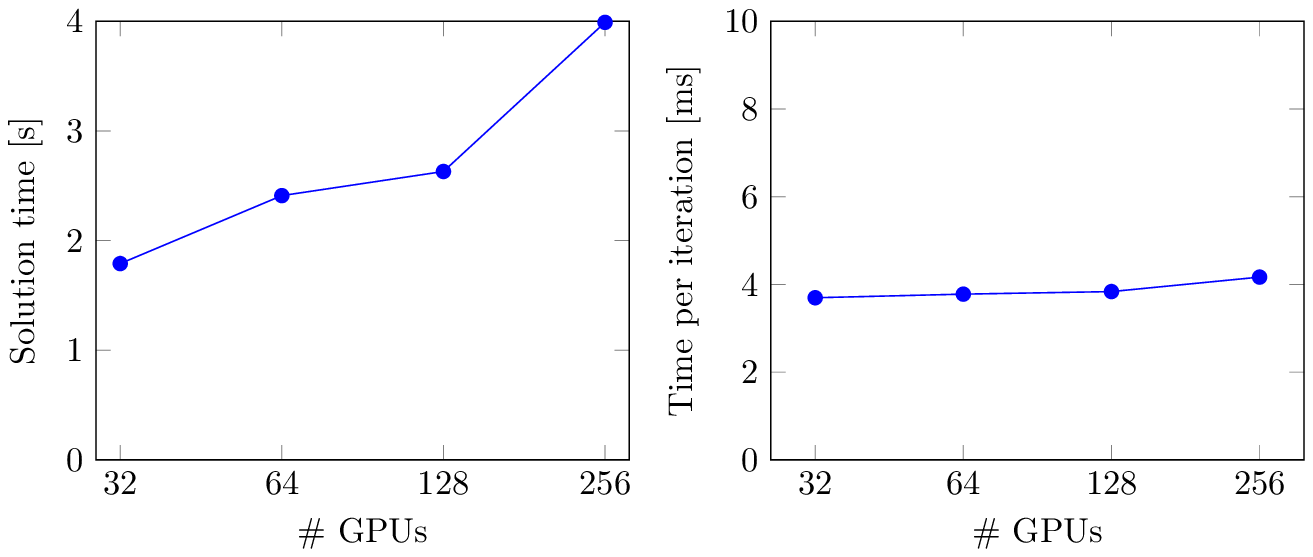,height=5cm}
\caption{Weak scalability of aFSAI CG on the Poisson model problem.
Number of CG iteration (top left), iteration time [s] (top right), set-up time [s]
(bottom left) and time for one iteration [s] (bottom right) vs. the number
of GPUs.}
\label{weakPois}
\end{center}
\end{figure}

Finally, we test the weak scalability of our implementation using a 7-point stencil
Finite Differences discretization of the Poisson problem on a cubic domain. In this
experiment we keep a constant number of 1,771,561 equations for each GPU and
increase the number of GPUs from 32 to 256. For larger problem sizes, the number
of CG iterations increases (see Figure~\ref{weakPois} upper-left), since aFSAI is not an
optimal preconditioner, as for instance AMG. As a consequence, the iteration time
as well as the total time correspondingly increase. However, if we focus on the set-up
stage and the time required by a single CG iteration (bottom-left and bottom right of
Figure~\ref{weakPois}, respectively) they do not change or only slightly
change when the problem size increases, showing an almost perfect weak
scalability.

\section{Conclusions}
We presented a multi-GPU implementation of the adaptive FSAI
preconditioner aimed at reducing the number of conjugate gradient
iterations required for the solution of large sparse linear systems
originated by the discretization of PDE in a number of scientific and
engineering fields.  We proposed a new memory layout for the matrix
that facilitates the overlap between communication and computation.
We carried out an extensive set of numerical tests that show the
efficiency of the proposed solution.  We make available the library
including both the adaptive FSAI preconditioner and the conjugate
gradient solver so that other people may use it both for benchmarking
and production purposes. The software is available along with the
source of a sample program from https://www.m3eweb.it/chronos/.
For the future we expect to develop a new version that will provide
support for other forms of communication among GPU, for instance
the third generation of Nvlink communication technology
recently announced by Nvidia.

\bibliographystyle{siam}
\bibliography{biblio}

\end{document}